\documentstyle[11pt]{article}

\title{{\bf Convexity of domains of  Riemannian manifolds }}

\author{Rossella Bartolo\thanks{Part of the Ph.D. thesis of this author is
contained in this paper}
\\Departamento de Geometr\'{\i}a y Topolog\'{\i}a\\
Fac. Ciencias, Univ. Granada\\ Avenida Fuentenueva s/n\\18071 Granada Spain\and
Anna Germinario\\Dipartimento Interuniversitario di Matematica\\
Universit\`a degli Studi di Bari\\Via E. Orabona, 4\\
70125 Bari Italy\and
Miguel S\'anchez\thanks{Partially supported by MEC Grant PB97-0784-C03-01}
\\Departamento de Geometr\'{\i}a y Topolog\'{\i}a\\
Fac. Ciencias, Univ. Granada\\ Avenida Fuentenueva s/n\\18071 Granada Spain}
\date{  }
\makeatletter

\@addtoreset{equation}{section}

\newtheorem{defn}{Definition}

\@addtoreset{defn}{section}

\newtheorem{teo}[defn]{Theorem}

\newtheorem{lem}[defn]{Lemma}
\newtheorem{rem}[defn]{Remark}
\newtheorem{rems}[defn]{Remarks}
\newtheorem{cor}[defn]{Corollary}
\newtheorem{prop}[defn]{Proposition}

\makeatother

\newcommand{\D}{{\cal D}}
\newcommand{\M}{{\cal M}}
\newcommand{\be}{\begin{equation}}
\newcommand{\ee}{\end{equation}}
\newcommand{\Om}{\Omega^1}
\newcommand{\R}{{\bf R}}
\newcommand{\N}{{\bf N}}
\newcommand{\dimo}{{\bf Proof: }}
\newcommand{\inte}{\int_{0}^{1}}
\newcommand{\eps}{\epsilon}
\newcommand{\<}{\langle}
\renewcommand{\>}{\rangle}
\renewcommand{\(}{\left(}
\renewcommand{\)}{\right)}

\newcommand{\xm}{x_m(s)}
\newcommand{\xms}{x_m(s_m)}

\newcommand{\yms}{y_m(s_m)}
\newcommand{\integ}{\int_{s_m}^{s}}

\begin{document}
\maketitle


\newpage

\section{Introduction}

In this paper we shall discuss the problem of the geodesic connectedness of subsets of Riemannian manifolds. In particular, we shall prove the geodesic connectedness of  open domains (i.e. connected open subsets) $\D$ of a smooth Riemannian manifold $ ( \M, \< \cdot,\cdot \> )$, under reasonable assumptions; moreover, in the more relevant cases $\D$ will be shown to be {\em convex}, i.e.   any pair of its points can be joined by a (non necessarily unique) minimizing geodesic\footnote{The word ``convex" is used in different non-equivalent ways in the literature. Sometimes, it is reserved for open domains such that each two points can be joined by an {\it unique} minimizing geodesic; for those following this convention, a better name for our domains would be {\it weakly convex}.}. As pointed out by Gordon \cite{gor}, this problem is important not only by its own but also because of its relation, via the Jacobi metric, to the problem of connecting two points by means of a trajectory
of
fixed energy for a Lagrangian system. Until now
this topic has been faced by using different techniques and under
different assumptions which allow us to control the non--completeness
of $\D$. In particular, geodesic connectedness can be
proved by using variational methods. In this case, the right assumption to get  existence, and multiplicity in same
cases, of geodesics connecting two fixed points, is a
convexity assumption on the boundary of $\D$ (see e.g. \cite{sal1});
here we shall work under weaker assumptions.

Our study makes necessary to discuss the different notions of
convexity for the (possibly singular) boundary points of the open
domain $\D$. In what follows, differentiability will mean
${\cal C}^4$; indeed, we shall need just ${\cal C}^3$ for
our main results (cfr. Theorems \ref{t1},  \ref{t2},
and ${\cal C}^2$ for Theorem \ref{t0}), but in our references the
highest assumption of differentiability is ${\cal C}^4$;  so, we prefer stating the results  assuming it. At first, recall that, by the well--known Hopf--Rinow  theorem, if a Riemannian manifold $\M$ is  complete then it is geodesically connected. From a variational point of view, the Hopf--Rinow theorem can be easily proved by using the functional
\be\label{0.1}
f(x) = \frac{1}{2} \inte\<\dot x(s), \dot x(s)\>ds
\ee
defined on a suitable Hilbert manifold, see Section \ref{2}.
It is well--known that the critical points of $f$ are geodesics and
it is not difficult to prove that $f$ admits a minimum point.
Then, in the complete case, this result guarantees not only
that the manifold is geodesically connected, but also that it is convex.

We examine now the case when the boundary $\partial \D$ of $\D$ in $\M$ is differentiable, that is,
$\overline \D = \D \cup \partial\D$ is a Riemannian manifold with (differentiable) boundary. Recall the following two natural notions of convexity around a point of the boundary.
\begin{defn}[Infinitesimal convexity]\label{d0.1}
We say that $\partial \D$ is infinitesimally  \linebreak convex at
$p\in \partial \D$ if the second fundamental form $\sigma_p$, with respect
to the interior normal, is positive semidefinite.
\end{defn}
\begin{defn}[Local convexity]\label{d0.2}
We say that $\partial\D$ is locally convex at $p\in\partial \D$ if there
exists a neighborhood $\overline U \subset \overline \D$ of $p$ such that
\be\label{0.2}
\exp_p\left( T_p\partial\D\right)\cap\left(\overline U\cap\D\right) = \emptyset.
\ee
\end{defn}
It is not difficult to show that the local convexity implies
the infinitesimal one, but the converse is not true.
Nevertheless, if the infinitesimal convexity is assumed on
a neighbourhood of a point of the boundary, then the notions are equivalent, as proved by Bishop \cite{bis}.
In order to apply variational methods to the study of geodesic
connectedness, a characterization of the infinitesimal convexity
is useful. Firstly, note that, by the differentiability of the boundary,  for each $p\in\partial\D$ there exist
a neighborhood $ U\subset \M$ of $p$
and a differentiable function
$\phi: U \cap \overline \D \longrightarrow\R$ such that
\be\label{0.3}
\left\{\begin{array}{lll}
          \phi^{-1}(0) =  U\cap\partial\D \\
          \phi > 0 &    \mbox{on $U\cap\D$ }\\
          \nabla\phi(q)\not= 0  & \mbox{for any $q\in U\cap\partial \D.$}
         \end{array}
\right.
\ee
Then, it is easy to check  that
  $\partial\D$ is infinitesimally convex  at $p\in\partial\D$ if and only if
for one (and then for all) function $\phi$ satisfying (\ref{0.3})
we have
\be\label{0.4}
H_\phi(p)[v,v]\leq 0 \quad \forall v\in T_p\partial\D.
\ee

Now, we shall go from considerations around a point of the boundary,
through considerations on all the boundary, and, so, the completeness of $\overline \D$ becomes essential.  Note that if $\M$ is complete then so is $\overline \D$ and, even though the converse is not true, there is no loss of generality assuming it, because the Riemannian metric can be modified out of $\overline \D$ to obtain completeness (see for example \cite{NO}). By  standard arguments the function $\phi$ in (\ref{0.3}) can be found on all $\overline \D$, and, thus, we have the following equivalent definitions for the convexity of all the boundary.
\begin{defn}[Global convexity, variational point of view]\label{cb.v}
Assume that $\M$ is complete.
$\partial\D$ is convex if and only if for
one, and then for all, nonnegative function $\phi$ on
$\overline\D$ such that
\be\label{0.5}
\left\{\begin{array}{lll}
            \phi^{-1}(0) = \partial\D \\
            \phi > 0 &  \mbox{on $\D$} \\
            \nabla\phi(q)\not= 0, & \mbox{for any $q\in\partial \D$}
         \end{array}
\right.
\ee
we have
\be\label{0.6}
H_\phi(q)[v,v]\leq 0 \quad \forall q\in\partial\D,
v\in T_q\partial\D.
\ee
\end{defn}
It is worth pointing out that condition
(\ref{0.6}) is equivalent to a geometric notion of convexity.
Indeed the following definition is equivalent too, see \cite{ger}.
\begin{defn}[Global convexity, geometrical point of view]\label{pv}
Assume  that $\M$ is complete.
$\partial\D$ is convex if for any $p,q\in\D$ the range of any
geodesic $\gamma:[0,1]\longrightarrow\overline\D$ such that
$\gamma(0) = p, \gamma(1) = q$  satisfies
\be\label{0.7}
\gamma\([0,1]\)\subset\D.
\ee
\end{defn}
Condition (\ref{0.7}) is a generalization to
Riemannian manifolds of the usual notion of convexity given
in Euclidean spaces. Moreover, all the above conditions provide different ways to prove that, when $\M$ is complete:
\begin{center}
$\D$ is convex if and only if $\partial \D$ is convex.
\end{center}
In fact, a variational technique based on the use of the functional
(\ref{0.1})  and a penalization argument make possible to prove
that if $\partial\D$ is convex, then $\D$ is convex. On the other hand, by
using Definition \ref{pv} it is easy to show that if $\partial \D$ is not convex then neither is $\D$.

Now we are ready to examine the general case where $\partial \D$ is not differentiable or $\overline \D$ is not complete. By using the results above, it is
clear that if there exists a sequence
\[
\(\overline \D_m\)_{m\in\N}
\]
of complete submanifolds with convex (differentiable) boundary such that
\begin{equation}
\label{seq}
\overline\D_m\subset \overline\D_{m+1} \quad \hbox{and} \quad
\D = \bigcup_{m\in\N} \overline\D_m,
\end{equation}
 then $\D$ is geodesically
connected. As a first question we can wonder if $\D$ must be convex. In Section \ref{s.g}
 we answer this question, by showing that if $\overline \D$ is complete then $\D$ is convex.
More precisely, let $\overline \D^c$ be the canonical completation of $\D$
by using Cauchy sequences, and $\partial_c \D$ the corresponding boundary points,
 $\overline \D^c = \D \cup \partial_c \D$
($\overline \D^c$ is always complete as a metric space, but the boundary points in
$\partial_c \D$ are not necessarily differentiable and, if they are,
the metric may be non--extendible or degenerate). Note that any point
of $\partial \D$  naturally determines one or more points in $\partial_c \D$, and $\overline \D$ is complete if and only if all the points in  $\partial_c \D$ are of this type (in this case, we can assume that $\M$ is complete). Then we will prove:

\begin{teo}\label{t0}
Assume that $\overline \phi: \overline \D^c \rightarrow [0, \infty]$ is a continuous function such that:
\begin{list}{(\roman{enumi})}{\usecounter{enumi} \labelwidth 3 em
\itemsep 0pt \parsep 0pt}
\item $\overline \phi^{-1}(0) = \partial_c \D$;
\item there exists an infinitesimal  sequence $(a_m)_{m\in\N}$ such that
 $\overline \phi^{-1}( ]a_m, \infty])$ is a Riemannian submanifold with convex (differentiable) boundary
$\overline \phi^{-1}(a_m)$ for any $m\in\N$.
\end{list}
Then $\D$ is geodesically connected.

Moreover, if $\M$ is complete then $\D$ is convex.
\end{teo}
We point out that in order to obtain convexity, the assumption of completeness on $\M$ cannot be removed, as we shall shown in Section \ref{s.g} by a counterexample.

This result is proved by using geometrical methods, and it
makes possible to  generalize the results by Gordon in  \cite{gor} by showing that his hypotheses imply those in Theorem \ref{t0} (Section \ref{s.g}).

Our main objective in this paper is achieved in Section 3, where we use variational methods (under the natural assumption of completeness for $\M$) to show that $\D$ is convex
when there exists a sequence $\(\D_m\)_{m\in\N}$ of open domains invading $\D$ as in (\ref{seq}), whose boundaries are differentiable but not necessarily  convex (and, thus, each $\D_m$ may be non--geodesically connected), if a suitable estimate of the loss of convexity of $\partial\D_m$ and boundness of the sequence is assured. More precisely, the following result will be proved.
\begin{teo}\label{t1}
Let $\M$ be a complete Riemannian manifold and $\D$ an open domain of $\M$.
Assume that
there exists a positive differentiable function
$\phi$ on $\D$  such that
\begin{list}{(\roman{enumi})}{\usecounter{enumi} \labelwidth 3 em
\itemsep 0pt \parsep 0pt}
\item
$\lim_{x \rightarrow \partial \D}
\phi(x) = 0$;
\item
each $y\in \partial \D$ admits a neighbourhood $U\subset \M$ and constants $a, b> 0$ such that
$$a \leq \| \nabla \phi(x)\| \leq b \quad \forall x\in\D \cap U; $$
\item
the first and second derivatives of the normalized flow of
$\nabla \phi$ are locally bounded close to $\partial \D$, that is:
each $y\in \partial \D$ admits a neighbourhood $U\subset \M$
such that the induced local flow on $\D \cap U$ has first and second derivatives with bounded norms;
\item
there exist a decreasing and infinitesimal  sequence
$(a_m)_{m\in \N}$ such that each $y\in \partial \D$ admits a neighbourhood $U\subset \M$ and a constant $M \in \R$ satisfying:
\be\label{1.1}
H_\phi(x)[v,v]
\leq M \<v,v\>\phi(x) \quad
\forall x\in \phi^{-1}(a_m) \cap U,
v\in T_x\phi^{-1}(a_m), m\in \N.
\ee
\end{list}
Then $\D$ is convex.
Moreover if $\D$ is not contractible in itself, then for any $p, q\in\D$ there exists a sequence $(x_m)_{m\in\N}$ of geodesics in $\D$ joining them such that
$$\lim_{m \rightarrow \infty} f(x_m ) = \infty.$$
\end{teo}
\begin{rems}
\label{rem1.7}
{\em When the boundary $\partial \D$ is smooth and convex in
the sense of Definition \ref{cb.v}, the function $\phi$
in (\ref{0.5}) always satisfies all conditions {\it (i)--(iv)} above. Let us examine the role of each one of these hypotheses.

(1)
Hypotheses {\it (i), (iv)}: they imply first that, taking  $\D_m = \phi^{-1}(]a_m,\infty[)$, condition (\ref{seq}) is satisfied, and, second, that  even when the boundaries $\partial \D_m$ may be non--convex, their loss of convexity is (locally) bounded by (\ref{1.1}).

(2) Hypothesis {\it (ii)}:
from Theorem \ref{t0}, if {\it (iv)} is satisfied with $M\leq 0$  then {\it (ii)} can be replaced just by:
$$\|\nabla \phi(x)\| \neq 0, \;
\forall x\in \phi^{-1}(a_m),\, \forall m\in \N.$$
But when $M>0$ the hypothesis {\it (ii)} must be imposed to make {\it (iv)} meaningful. The reason is that the left hand side in (\ref{1.1}) describes the shape of $\phi^{-1}(a_m)$, but the value of $\phi$ in the right hand side can be almost arbitrarily changed, if no bound on the gradient is imposed.
More precisely, consider any  smooth function $\varphi: ]0,\infty[
\rightarrow ]0,\infty[$ such that $\lim_{s\rightarrow 0} \varphi(s) =0$ and its derivative satisfies $\dot \varphi>0$. Then $\phi^* = \varphi \circ \phi$ also satisfies {\it (i)}, and
$$
H_{\phi^*}(x)[v,v] = \dot \varphi(\phi(x)) H_\phi (x) [v,v]\quad
\forall x\in \phi^{*-1}(a_m^*),
v\in T_x\phi^{*-1}(a_m^*),
$$
where $a_m^* = \varphi (a_m), m\in \N$. When $\phi$ satisfies {\it (ii)} then
$\phi^*$ satisfies {\it (ii)} if and only if
$a^* \leq \dot \varphi \leq b^*$, close to 0,  for some $a^*, b^* >0$.
In this case, $\phi$ satisfies {\it (iv)} if and only if so does $\phi^*$.
But if {\it (ii)} were not imposed, it would be possible that one of the
functions satisfies {\it (iv)} and the other does not\footnote{In fact, assume
that $\phi$ satisfies {\it (iv)}, but there exists an infinitesimal
sequence $(a_m)_{m\in \N}$ such that $H_\phi(x_m)[v_m,v_m] >0$ for some $x_m$
at each $\phi^{-1}(a_m)$ and unitary vector
$v_m\in T_{x_m}\phi^{-1}(a_m)$. Then take $(k_m)_{m\in \N}$
such that
$k_m > 0,$ $\(k_mH_\phi(x_m)[v_m,v_m]\)_{m\in \N}$ is not infinitesimal.
Clearly any $\varphi$ as above such that
$\dot \varphi (a_m) = k_m$ yields the required example.}. This shows that {\it (iv)} is not reasonable by itself as a measure of the loss of convexity of the hypersurfaces $\phi^{-1}(a_m)=\phi^{*-1}(a_m^*)$, being hypothesis {\it (ii)} natural.

{\sloppy (3) Hypothesis {\it (iii)}: the bounds on the normalized flow
(i.e. the flow of $\nabla \phi / \|\nabla \phi\|^2 $) are technical, and they express the unique control we impose on the intermediate hypersurfaces between two consecutive $\phi^{-1}(a_m)$. Note that if {\it (iii)} were not imposed then   hypersurfaces arbitrarily close to $\partial \D$ ``very distorted" by the flow could exist.

}}
\end{rems}

Technical condition {\it  (iii)} and even the completeness of the ambient manifold $\M$ can be weakened if {\it (iv)} is imposed on all points and directions enough close to the boundary. So, a straightforward consequence of the technique in the proof of Theorem \ref{t1} is the following result (compare with \cite{sal1}):

\begin{teo}\label{t2}
Let $\M$ be a Riemannian manifold, $\D \subset \M$ an open domain,
and $\overline \D^c = \D \cup \partial _c\D$ its canonical Cauchy completation. Assume that
there exists a positive differentiable function
$\phi$ on $\D$  such that
\begin{list}{(\roman{enumi})}{\usecounter{enumi} \labelwidth 3 em
\itemsep 0pt \parsep 0pt}
\item
$\lim_{x \rightarrow \partial_c \D}
\phi(x) = 0$;
\item
each $y\in \partial_c \D$ admits a neighbourhood $U\subset \overline \D^c$ and constants $a, b > 0$ such that
$$a\leq \|\nabla \phi(x)\| \leq b \quad \forall x\in U\cap \D;$$
\item
each $y\in \partial_c \D$ admits a neighbourhood $U\subset  \overline \D^c$ and
a constant $M\in \R$ such that inequality (\ref{1.1}) holds for all $x \in\D \cap U$ and for all $v \in T_x \M$.
\end{list}
Then $\D$ is convex.
Moreover if $\D$ is not contractible in itself, for any $p, q\in\D$ there exists a sequence $(x_m)_{m\in\N}$ of geodesics in $\D$ joining them such that
$$\lim_{m \rightarrow \infty} f(x_m ) = \infty.$$
\end{teo}

\begin{rems}
{\em Note that when $\partial_c \D$ is convex in the sense of
Definition \ref{cb.v} the function $\phi$ in (\ref{0.5})
does not necessarily satisfy {\it (iii)} because this condition is now
imposed on all tangent vectors $v$. Nevertheless,
when Definition \ref{cb.v} is applicable, it is independent of
the chosen $\phi$; thus, varying $\phi$, all tangent $v$ can be considered as tangent to a level hypersurface. So, we can conclude that the hypotheses in Theorems \ref{t1} and \ref{t2} imply, for each point of the boundary, local conditions which extend those for differentiable boundaries. }
\end{rems}

In Section \ref{2} Theorem \ref{t1} will be proved. To this aim, we shall
penalize the
functional $f$ of (\ref{0.1}) with a term depending on a positive parameter
$\eps$ and we shall study the Euler--Lagrange equation associated to the penalized
functionals $f_\eps$. The crucial point is to prove that a critical point
of $f_\eps$ in a
sublevel of $f_\eps$ is uniformly far (with respect to $\eps$) from $\partial \D$.
In the proof we shall  ``project'' the critical points of the penalized
functionals (using the normalized flow of $\nabla \phi$) on the hypersurface $\phi^{-1} (a_m )$ for $m$ large enough. This makes possible to get critical points of $f$ (i.e. geodesics) not touching  $\partial \D$ by means of a limit process.

Finally, in Section \ref{Appendix} the discussion of Theorems \ref{t1}, \ref{t2}  is completed
by  giving: (a) some examples which show the applicability and independence of the
hypotheses of Theorems \ref{t0}, \ref{t1}, \ref{t2}, and (b) an application to the  existence of trajectories of fixed energy for dynamical systems.

\section{Proof of Theorem 1.5 and Gordon's theorem}
\label{s.g}

{\bf Proof of Theorem \ref{t0}.} Given $p, q \in \D$ choose
$a_m < \min\,\{ \overline \phi(p), \overline \phi(q) \}$. By {\it (i)} and
the continuity of
$\overline\phi$ at $\partial_c\D$, the manifold $\overline\D_m =
\overline\phi^{-1} ([a_m,\infty])$ is complete. Then  $p$ and $q$
belong to the interior
of $\overline\phi^{-1}([a_m, \infty])$, which is intrinsically convex.

Now let us show that if $\M$ is complete then $\D$ is convex.
Let $m_0\in\N$ be the first positive integer such that
$a_{m_0}<\min\{\overline\phi(p), \overline\phi(q)\}$.
As before we get for any $m\geq m_0$ the
existence of a minimizing geodesic $\gamma_m$ in $\D_m$ joining $p$ and $q$.
Let us consider $\gamma_m$ parametrized by arc length, so
$\gamma_m:[0,l_m]\longrightarrow \D_m$,
$\<\dot\gamma_m(s),\dot\gamma_m(s)\> = 1$, for any $s$.
Note that, up to a subsequence,
$\dot\gamma_m(0)\longrightarrow v$, with $v$ unitary vector, and that
$(l_m)_{m\in\N}$ is a decreasing sequence converging to the distance $l$ between $p$ and $q$.
By standard arguments, the geodesic
$\gamma:[0,l]\longrightarrow \M$ with $\dot\gamma(0) = v$ has range in $\D\cup\partial\D$, joins $p$ and $q$, and satisfies
$$\gamma_m\longrightarrow \gamma\quad\hbox{uniformly in $[0,l]$.}$$
If, for infinitely many $m\in\N$, $\gamma\equiv \gamma_m$ the proof is complete.
So let us assume $\dot\gamma_m(0)\not= v$
for infinitely many $m$ and $(l_m)_{m\in\N}$ strictly decreasing: note that,
necessarily, $p$ and $q$ are conjugate along $\gamma$.
Let $U$ be a star--shaped neighbourhood of $q$, with $U\subset\D_{m_0}$.
We shall prove that, for small $\delta$,  $q_\delta= \gamma(l-\delta)\in U$
is always conjugate
to $p$ along $\gamma$ getting a contradiction.
Since $p,q_\delta\in \D_m$ for any $m\geq m_0$, reasoning as before, we can
find a sequence of minimizing geodesics parametrized by arc length
$\bar\gamma_m:[0,\bar{l}_m]\longrightarrow \D_m$
joining $p$ and $q_\delta$,
a subsequence
$\dot{\bar\gamma}_m(0)\longrightarrow \bar v$, with $\bar v$ unitary
vector (which we can assume distinct from $v$, otherwise $q_\delta$
 is also conjugate and the proof is complete) and
$\bar\gamma:[0,\bar l]\longrightarrow \D\cup\partial\D$
with $\dot{\bar\gamma}(0) = \bar v$.
Again
$(\bar l_m)_{m\in\N}$ is decreasing and converging to the distance
$\bar l$ between $p$ and $q_\delta$. We claim $\bar l \leq l - \delta$;
otherwise, modifying slightly the curves $\gamma_m$, we could find a curve
joining $p$ and $q_\delta$ contained in $\D_m$ for $m$ large, with length less than $\bar l$, which is a contradiction. Now we can define the union curve
of $\bar\gamma$
on $[0,\bar l]$ and
$\gamma$ restricted to $[l - \delta, l]$
which
joins $p$ and $q$ and with length less or equal than $l$.
As this union curve is not differentiable at $\bar l$, we can slightly
modify it to obtain a curve $\hat\gamma$ with length less than $l$ and which agrees with $\bar\gamma$ out of $U$.
Finally, we could find for large $m$ a curve joining $p$ and $q$, with length less than $l$ and which agrees with some $\bar\gamma_m$ out of $U$, getting an absurd with the minimality of $\gamma_m $ since $U\subset \D_m$.
$\Box$

\vspace{0.5 cm}
At the end of this section we shall give a counterexample for the case $\M$ non complete.

\vspace{0.5 cm}
The possibility of extending Gordon's results by using variational methods has already been pointed out in \cite[Chapter 4]{ma}, \cite{sal1}. Nevertheless our point of view is quite different, and the extension we obtain is, at any case, elementary and stronger. In fact,  in Gordon's result a convexity assumption is done on the whole manifold; in the quoted references it is claimed that this global assumption must imply a convexity property close to the boundary, which should be enough from a variational point of view; finally, we will check now that the global assumption imply a convexity property for a sequence of hypersurfaces close to the boundary, which is enough from any of the points of view sketched in Section 1.

We recall that a map $h$ between manifolds is said to be {\it proper} if
$h^{-1}(K)$ is compact whenever $K$ is compact. In particular, if
$h:\D \rightarrow \R$ is proper necessarily $|h(p)|\rightarrow \infty $ as
$p \rightarrow \partial \D$. Recall also that a (real--valued ${\cal C}^2$)
function is called {\it convex} when its Hessian is positive semidefinite. Gordon's result \cite[Theorem 1]{gor} asserts:

\begin{teo}
\label{gordon}
If the open domain $\D$ of the Riemannian manifold $\M$ supports
a proper positive convex function $h$, then it is geodesically connected.
\end{teo}

We will reprove this result, by showing that its hypotheses imply the ones
in  Theorem \ref{t0}. There is a second theorem in Gordon's paper, which
can be reproved in the same way. Recall first that Theorem \ref{t0} can also
be stated  assuming $\overline\phi^{-1}(\infty ) = \partial_c \D$
and $(a_m)_{m\in\N}$ diverging to $\infty$.

\vspace{0.5 cm}
\noindent {\bf Proof of Theorem \ref{gordon}.} By Sard's theorem, almost all
the values of $h$ are regular. Thus, as $h$ is proper (and positive),
there exists a diverging sequence of regular values $(a_m )_{m\in\N}$
contained in the range of $h$. Moreover,
$h^{-1} ([0,a_m])$ is a compact Riemannian manifold with boundary $h^{-1}(a_m)$. From Definition \ref{cb.v}, this boundary is convex (put $\phi = h(a_m)-h$). So,
extending continuously $h$ to a function
$\overline h : \overline \D^c \rightarrow [0,\infty]$,
Theorem \ref{t0} (in the version above) can be claimed. $\Box$

\vspace*{0.5 cm}
\noindent {\bf A counterexample.} The following counterexample shows that the result in Theorem \ref{t0} on geodesic connectedness cannot be strengthened to obtain convexity, when $\M$ is not complete.

Consider two open hemispheres $H_0, H_1$ in $\R^3$ and let $x_0 , x_1$ be
their north poles (see Fig. 1). Put a sequence of inmersed tubes
$(T_m)_{m\in\N}$ connecting $H_0$ and $H_1$ of decreasing length and
such that any curve joining $x_0$ and $x_1$ through $T_m$ is longer
than a minimizing curve joining them through $T_{m+1}$. We also assume that
the width of these tubes goes to zero, and   their mouths in each
hemisphere go to a point $e_i, i=0,1$ in the equator, being all their
centers in the same meridian (the shape of the resulting hemispheres is
shown in Fig. 2). Let $\M$ be this manifold inmersed in $\R^3$, and $\D = \M$.

Recall that $\partial_c \D$ is canonically identifiable to the equators, and let $\phi : \D \rightarrow \R$ be the height function with
$\lim_{p\rightarrow \partial_c \D}\phi(p)=0$. Clearly, a sequence $(a_m)\rightarrow 0$ can be chosen such that $\overline\D_m\ = \phi^{-1}([a_m, \infty[)$ is a complete Riemannian manifold with convex boundary $\partial\D_m = \phi^{-1}(a_m)$,  containing the tube $T_m$ but not $T_{m+1}$, and satisfying (\ref{seq}).
So, in each convex manifold $\D_m = \phi^{-1}(]a_m, \infty[)$ there exists a minimizing geodesic $\gamma_m$ connecting $x_0 , x_1$, and, by the condition on the lengths of the tubes,
\[
\hbox{length}\,(\gamma_{m+1}) < \hbox{length}\,(\gamma_m).
\]

So, if there was a minimizing geodesic $\gamma$ between $x_0 , x_1$ in $\D$, necessarily it should be included in some $\D_m$ and
$\hbox{length}\,(\gamma_{m+1}) < \hbox{length}\,(\gamma)$,
a contradiction.

\section{Proof of Theorems 1.6 and 1.8 } \label{2}
Before introducing the functional framework, we recall that,
by the well--known
Nash embedding Theorem (see \cite{na}),
any smooth Riemannian manifold  $\M$ is isometric to a
submanifold of ${\bf R}^N$, with $N$ sufficiently large, equipped with the
metric induced by
the Euclidean metric in ${\bf R}^N$. So, henceforth, we shall assume
that $\M$ is a submanifold of ${\bf R}^N$ and $\< \cdot,\cdot \>$ is the Euclidean metric.
It is well--known that      the
geodesics in $\D$ joining two fixed   points $p$ and $q$ of $\D$  are the critical points of the
{\em action integral}  (\ref{0.1})
defined on $\Om\(\D\)$
where
$$\Om\( \D\) = \left\{ x\in H^{1,2} ([0,1], \D) \mid
x(0)= p, x(1) = q\right\}$$
and
$$H^{1,2} ([0,1], \D) = \left\{ x\in H^{1,2} ([0,1], \R^N ) \mid
x([0,1]) \subset \D \right\}.$$
It can be proved that $\Om\(\D\)$ is a Hilbert submanifold of
 $H^{1,2} ([0,1], \D)$ whose tangent space at $x\in \Om\(\D\)$
  is given by
$$T_x\Om\(\D\)= \left\{ v \in H^{1,2}([0,1], T\D) \mid
v(s) \in T_{x(s)}\D \;\; v(0) = 0 = v(1)\right\}.$$

\vspace{3mm}
We recall the following definition.
\begin{defn}\label{ps}
Let $(X,g)$ be a Riemannian manifold
modelled on a Hilbert space and let $F\in {\cal C}^1(X,{\bf R})$. We say
that $F$ satisfies the {\em Palais--Smale condition} if every
sequence $(x_m)_{m\in \N}$ such that
\be \label{2.2}
(F(x_m))_{m \in {\bf N}} \hbox{ is bounded},
\ee
\be \label{2.3}
\|\nabla F(x_m)\| \rightarrow 0,
\ee
contains a converging subsequence, where $\nabla F(x)$ denotes the gradient
of $F$ at the point $x$ with respect to the metric $g$
and $\Vert \cdot \Vert$ is the norm on the tangent bundle induced by $g$.
A sequence satisfying (\ref{2.2})--(\ref {2.3}) is said a Palais--Smale
 sequence.
\end{defn}
In our case there are Palais--Smale sequences
that could converge to a curve which ``touches" the boundary
$\partial \D$, so we penalize the functional $f$ in a suitable way,
following \cite{bfg1}.
For any $\eps\in ]0,1]$, we consider on $\Om\(\D\)$
the functional
\be \label{feps}
 f_\eps(x) =
f(x) + \inte{\eps\over {\phi^2(x)}}ds
\ee
where $\phi$ has been introduced in Theorem \ref{t1}.
For any $\eps\in ]0,1]$ $f_\eps$ is a ${\cal C}^2$ functional
and if $x\in\Om\(\D\)$ is a critical point of $f_\eps$, by using a
boot--strap argument it can be proved that it is ${\cal C}^2$
and satisfies
\be\label{2.5}
D_s\dot x =
- {{2\eps}\over{\phi^3(x)}}\nabla\phi(x).
\ee
For any $\eps\in ]0,1]$, $s\in [0,1]$, we set
\be\label{3.1}
\lambda_\eps(s) =
{{2\eps}\over{\phi^3(x(s))}},
\ee
which represents the multiplier in (\ref{2.5}).
Multiplying (\ref{2.5}) by $\dot x$, it is easy to get the existence
of a constant $E_\eps(x)\in\R$ such that
\be\label{2.6}
{1\over 2}\<\dot x(s),\dot x(s)\> -
{{\eps}\over{\phi^2(x(s))}} =  E_\eps(x) \quad \forall s\in [0,1].
\ee
To prove that the penalized functionals satisfy the Palais--Smale
condition, we recall the following lemma which holds with slight variants of the proof in \cite{bfg1} also under our local assumptions {\it (i)--(ii)} of Theorem \ref{t1}.
\begin{lem}\label{l1}
Let $(x_m)_{m\in \N}$ be a sequence in $\Om\(\D\)$
such that
\be \label{2.7}
\sup_{m\in \N}\inte \< \dot x_m , \dot x_m \> ds < \infty
\ee
and assume the existence of a sequence $(s_m)_{m\in \N}$ in $[0,1]$ such that
\be \label{2.8}
\lim_{m\longrightarrow \infty} \phi(x_m(s_m)) = 0.
\ee
Then
\be \label{2.9}
\lim_{m\longrightarrow  \infty}
\inte{1\over {\phi^2(x_m(s))}}~ds = \infty.
\ee
\end{lem}
\begin{prop}\label{p1}
Let $f_\eps$ be as in (\ref{feps}). Then
\begin{list}{(\roman{enumi})}{\usecounter{enumi}\labelwidth  5em
\itemsep 0pt \parsep 0pt}
\item
for any
$\eps \in]0,1]$ and for any $c \in {\bf R}$ the sublevels
$$
f_{\eps}^c = \{ x \in \Om\(\D\) | f_\eps(x) \leq c \}
$$
are complete metric subspaces of $\Om\(\D\)$;
\item
for any $\epsilon \in ]0,1]$, $f_\epsilon$ satisfies the Palais--Smale
condition.
\end{list}
\end{prop}
\dimo
For any  $\eps\in ]0,1], c\in\R$, let $(x_m)_{m\in \N}$ be a Cauchy sequence
in $f_\eps^c$, then
it is a Cauchy sequence also in
$H^{1,2}([0,1],{\bf R}^N)$, so it converges strongly
to a curve $x$ in $H^{1,2}([0,1],{\bf R}^N)$.
Since this convergence is also uniform, by Lemma \ref{l1} it results that
$x \in \Om\(\D\)$ and by the continuity of $f_\eps$, we obtain
the
first part of the proposition.
Now let $(x_m)_{m\in \N}$ be a Palais--Smale sequence; in particular it results
that
\be
\inte \< \dot x_m , \dot x_m \> ds
\quad\hbox{is bounded}.
\ee
Then, up to a subsequence, we
get the existence of a $x \in H^{1,2}([0,1],{\bf R}^N)$
such that
\be\label{2.10}
x_m \longrightarrow x \hbox{ weakly in } H^{1,2}([0,1],{\bf R}^N).
\ee
Arguing as in the first part of the proof, we get that
$x \in \Om\(\D\)$.
Using standard arguments, it can be proved that
$$
x_m\longrightarrow x  \hbox{ strongly in } H^{1,2}([0,1],{\bf R}^N).\quad\Box
$$
\begin{rem}\label{r1}
{\em By Proposition \ref{p1}, for any $\eps\in]0,1]$, $f_\eps$ has a minimum
point $x_\eps\in\Om\(\D\)$; it is easy to see that there exists $k>0$ such
that
$$
f_\eps(x_\eps)\leq k,
$$
for any $\eps\in]0,1]$. Moreover, by (\ref{2.6}) we get for any $\eps\in]0,1]$
$$
E_\eps(x_\eps) =
f_\eps(x_\eps) - 2\inte{\eps\over {\phi^2(x_\eps(s))}}~ds \leq k,
$$
hence
\be\label{2.11}
{1\over 2}\<\dot x_\eps(s),\dot x_\eps(s)\>
\leq k + {\eps\over {\phi^2(x_\eps(s))}},
\ee
for any $\eps\in]0,1]$, $s\in [0,1]$.}
\end{rem}
\begin{rem}\label{r3}
{\em
In the sequel, we shall need to relate the Hessian of a
$\hat\Psi\in {\cal C}^2(\R^N,\R)$ to the one of its
restriction $\Psi$ on $\M$.
For any $y\in \M$ let
\be\label{3.2}
P(y): {\R}^N\longrightarrow T_y {\cal M},
\ee
\be\label{3.3}
Q(y): {\bf R}^N\longrightarrow T_y  \M^\perp,
\ee
be respectively the projections on $T_y {\cal M}$ and
$T_y {\cal M}^\perp$.
Since
$ \M $ is a ${\cal C}^3$ submanifold of ${\bf R}^N$,  there  exist
$A^i_j\in {\cal C}^2( {\cal M},{\bf R})$,
$i$,$j \in \{1, \dots, N\}$ such that
for any $y\in \M, v=(v^1,\ldots ,v^N) \in {\bf R}^N$
$$
Q(y)[v] = \sum_{i,j=1}^N A^i_j(y)v^j e_i ,
$$
 where $e_1,..., e_N $ is the canonical basis of $\R^{N}$.
We locally extend the functions $A^i_j$ to ${\cal C}^2$ functions
(still denoted by $A^i_j$) on ${\bf R}^N$. For any $y \in  {\bf  R}^N$,  we
define the differential map
${\rm d}Q(y)\colon {\bf R}^N\times {\bf R}^N \longrightarrow{\bf R}^N$ as
$$
{\rm d}Q(y)[v,w] =
\sum_{i,j,k=1}^N {{\partial A_j^i(y)}\over\partial {x_k}}v^k w^j e_i,
\quad \forall v,w \in {\bf R}^N.$$
Even if
${\rm d}Q$ could depend on
the extensions of the functions $A^i_j$, for any
$y\in {\cal M}$,
 the restriction of ${\rm d}Q(y)$ to
$T_y  {\cal M}\times T_y  {\cal M}$ is well--defined.
It can be proved (see e.g. \cite[Lemma 8]{ger}) that for any
$y\in  {\cal M}$, $v\in T_y {\cal M}$:
$$
H_{\Psi}(y)[v,v] =
{\rm d}^2\hat\Psi(y)[v,v] -
{\rm d}\hat\Psi(y)[{\rm d}Q(y)[v,v]],
$$
where $d$ and $d^2$ are the
differential map and the second differential map on ${\bf R}^N$.}
\end{rem}
\begin{lem}\label{l3}
Let $\(x_\eps\)_{\eps > 0}$ be a family in $\Om(\D)$ of critical
points of $f_\eps$ such that
\be\label{2.12}
f_\eps(x_\eps)\leq k \quad\forall\eps\in ]0,1],
\ee
for a suitable positive constant $k$.
Then ${\displaystyle \(\lambda_\eps(s) = {{2\eps}\over{\phi^3(x_\eps(s))}}\)_{\eps >0}}$ is bounded in
$L^\infty\([0,1],\R\)$.
\end{lem}
\dimo
Let $\(\eps_m\)_{m\in\N}$ be a decreasing and infinitesimal  sequence
in $]0,1]$ and let
$\(x_{\eps_m}\)_{m\in \N}$ be a sequence of critical points of
$f_{\eps_m}$ satisfying (\ref{2.12}).
For the sake of simplicity, in the following we set
$x_{\eps_m} \equiv x_m$, $\lambda_{\eps_m} \equiv \lambda_m$.
Now let $u_m(s) = \phi(x_m(s))$, for $s\in [0,1], m\in\N$ and
$u_m(s_m) = \min_{s\in [0,1]}u_m(s)$, for any $m\in\N$.
It suffices to prove the lemma when,
up to a subsequence,
\be\label{3.4}
\lim_{m\longrightarrow \infty}u_m(s_m) = 0.
\ee
By (\ref{2.12}) and the Poincar\'e inequality there exists
$x\in H^{1,2}([0,1],{\bf R}^N)$ such that
\be\label{uni}
x_m\longrightarrow x \quad \hbox{uniformly},
\ee
and since, up to a subsequence,
\be\label{2.15}
\lim_{m\longrightarrow \infty}s_m = s_0,
\ee
by (\ref{3.4}) and (\ref{uni}) it easily follows
$s_0\in ]0,1[$ since $\phi(p), \phi(q) > 0$.
It is not difficult to see that $\(x_m (s_m)\)_{m\in \N}$ converges to $x(s_0)=y \in \partial \D$. Let $U$ be a neighborhood of $y$ such that {\it  (ii)}--{\it  (iv)} of Theorem \ref{t1} hold, then there exists $\mu > 0$ such that
\be \label{3.9}
x_m (s) \in U \cap \D
\ee
 for any $s\in J = [s_0 - \mu, s_0 + \mu]$ and $m$ sufficiently large.
By  (\ref{2.5}), (\ref{3.1}), (\ref{3.4}) and {\it  (ii)} of Theorem \ref{t1},
we get for m large enough
\begin{eqnarray} \label{3.5}
& \ddot u_m(s_m) =
H_\phi(\xms)[\dot \xms,\dot \xms] -
\lambda_m(s_m){\|\nabla\phi(\xms)\|}^2\leq & \nonumber \\
& H_\phi(\xms)[\dot \xms,\dot \xms]
- a^2\lambda_m(s_m). &
\end{eqnarray}
By the assumptions of Theorem \ref{t1} and (\ref{3.4}),
for any $m\in\N$ there exists $k_m\in\N$ such that
$$a_{k_m}\leq \phi(x_m(s_m))\leq \phi(x_m(s))$$
 for any $s\in [0,1]$.
Let us consider the Cauchy problem
\be\label{2.16}
\left\{
\begin{array}{ll}
{\displaystyle  \dot\eta = -\frac{ \nabla \phi(\eta)}{ \| \nabla
\phi(\eta ) \|^2}  }\\                                                         
\eta(0) = x \in   U\cap  \D                                                 
\end{array}
\right.
\ee
and call $\eta(s,x): {\cal U}\subset \R \times \overline\D
\longrightarrow \overline \D$ the flow associated to the Cauchy problem
(\ref{2.16}), where ${\cal U}$ is the maximal domain where the flow can be defined. Set for any $s\in J$
\begin{eqnarray}
\label{2.17}
\tau_m(s) = \phi(x_m(s)) - a_{k_m}\\
\label{2.18}
y_m(s) = \eta\(\tau_m(s), x_m(s)\).
\end{eqnarray}
Observe that if $m$ is sufficiently large and $\mu$ is opportunely chosen
$$( \tau_m (s) , x_m(s) ) \in {\cal U}\qquad  y_m (s) \in U \qquad \forall s\in J.$$
Then, for $m$ large enough, we can define the projection
$\Pi_m : U\cap \D  \longrightarrow\phi^{-1}(a_{k_m})$
\be\label{2.19}
\Pi_m(x_m(s)) = y_m(s)\quad s\in J.
\ee
Note that, by the definition of $\Pi_m$,
$\<\dot y_m(s),\nabla\phi(y_m(s))\> = 0$, for any
$s\in J$, so since
$$
\dot y_m(s) =
\eta_x\(\tau_m(s), x_m(s)\)[\dot x_m(s)] -
{{\nabla\phi(y_m(s))}\over {{\|\nabla\phi(y_m(s))\|}^2}}
\dot u_m(s),
$$
by assumption {\it (iii)} of Theorem \ref{t1} we get for any $s\in J$
\begin{eqnarray} \label{2.27}
& {\displaystyle
{|\dot y_m(s)|}^2 =
{| \eta_x\(\tau_m(s), x_m(s)\)[\dot x_m(s)]|}^2
-
{{\dot u_m^2(s)}\over {{\|\nabla\phi(y_m(s))\|}^2}}\leq } & \nonumber \\
 & {\displaystyle
C|\dot x_m (s) |^2 .} &
\end{eqnarray}
Hence by
(\ref{3.5}) and {\it  (iv)} of Theorem \ref{t1}
 \begin{eqnarray}  \label{2.20}
& {\displaystyle
a^2\lambda_m(s_m)\leq
H_\phi(\xms)[\dot \xms,\dot \xms] -
H_\phi(\yms)[\dot \yms),\dot \yms] + } & \nonumber \\
& {\displaystyle
M\<\dot \yms,\dot \yms\>a_{k_m}. }&
\end{eqnarray}
Now, following Remark \ref{r3}, set for any
$x\in\D$, $v, w \in {\bf R}^N$
\be\label{3.6}
L(x)[v,w] =
{\rm d}^2\phi(x)[v,w] - {\rm d}\phi(x)[{\rm d}Q(x)[v,w]],
\ee
\begin{eqnarray}  \label{2.23}
& {\displaystyle
H_\phi(\xms)[\dot \xms,\dot \xms] -
H_\phi(\yms)[\dot \yms,\dot \yms] = } & \nonumber \\
& {\displaystyle
L(\xms)[\dot \xms,\dot \xms] -
L(\yms)[\dot \yms,\dot \yms] =} & \nonumber \\
& {\displaystyle
L(\xms)[\dot \xms,\dot \xms] -
L(\yms)[\dot \xms,\dot \xms]+} & \nonumber \\
& {\displaystyle
L(\yms)[\dot \xms,\dot \xms] -
L(\yms)[\dot \yms,\dot \yms] =} & \nonumber \\
& {\displaystyle
L(\xms)[\dot \xms,\dot \xms] -
L(\yms)[\dot \xms,\dot \xms]+} & \nonumber \\
& {\displaystyle
L(\yms)[\dot \xms + \dot\yms,\dot \xms -\dot\yms].}
\end{eqnarray}
In the following we shall denote by $M_1, \ldots , M_7 $
suitable positive constants.
Since $\phi\in{\cal C}^3$, $\eta\in{\cal C}^2$, using the
mean value theorem, the boundness of
$\(\|x_m \|_\infty\)_{m\in \N}$ and $\(\|y_m \|_\infty\)_{m\in \N}$, {\it (iii)} of Theorem \ref{t1} and (\ref{3.4}):
\begin{eqnarray}  \label{3.8}
& {\displaystyle
L(\xms)[\dot \xms,\dot \xms] -
L(\yms)[\dot \xms,\dot \xms] \leq} & \nonumber \\
& {\displaystyle
M_1{|\dot\xms|}^2|\xms - \yms| \leq} & \nonumber \\
& {\displaystyle
M_2{|\dot\xms|}^2\(u_m(s_m) - a_{k_m}\)\leq} & \nonumber \\
& {\displaystyle
M_3{|\dot\xms|}^2.} &
\end{eqnarray}
By the boundness of $\(\|y_m \|_\infty\)_{m\in \N}$ and (\ref{2.27}) we get
\begin{eqnarray} \label{3.7}
& {\displaystyle
L(\yms)[\dot \xms + \dot\yms,\dot \xms -\dot\yms]\leq }& \nonumber \\
& {\displaystyle
M_4{|\dot\xms|}^2}. &
\end{eqnarray}
Then by
(\ref{2.20}), (\ref{2.23}), (\ref{3.8}),
(\ref{3.7}), (\ref{2.27}), we get
$$
a^2\lambda_m(s_m)\leq
M_5{|\dot\xms|}^2
$$
so that (\ref{2.11}) implies
$$
{\eps_m\over {\phi^3(x_m(s_m))}}\leq
M_6 + M_7{\eps_m\over {\phi^2(x_m(s_m))}}
$$
from which the boundedness of the multiplier follows. $\Box$
\begin{prop}\label{l2}
Let
$\(x_\eps\)_{\eps >0}$ be a sequence of critical points of
$f_\eps$ and $k$ a positive constant such that
(\ref{2.12}) holds.
Then there exists a positive constant $\beta$ such that
\be\label{2.13}
\phi(x_\eps(s))\geq \beta\quad \forall s\in [0,1], \eps\in]0,1].
\ee
\end{prop}
\dimo
Assume by contradiction that there exist a
decreasing and infinitesimal sequence $\(\eps_m\)_{m\in\N}\in ]0,1]$ and a sequence
$\(x_{\eps_m}\)_{m\in \N}$ of critical points of
$f_{\eps_m}$ satisfying (\ref{2.12}) and
such that
\be \label{2.14}
\min_{s\in [0,1]} \phi(x_{\eps_m}(s))\rightarrow 0 \hbox{ as }  m\longrightarrow \infty.
\ee
Defining $x_m, u_m, s_m, s_0, U, J$ and
$a_{k_m}$ as in Lemma \ref{l3},
since for $s\in J$, $m\in \N$ sufficiently large $y_m(s)\in U$,
we get, by  {\it  (iv)} of Theorem \ref{t1}
\begin{eqnarray}  \label{3.10}
& {\displaystyle
\ddot u_m(s)\leq
H_\phi(x_m(s))[\dot x_m(s),\dot x_m(s)] -
H_\phi(y_m(s))[\dot y_m(s),\dot y_m(s)] + } & \nonumber \\
& {\displaystyle
M\phi(y_m(s))\<\dot y_m(s),\dot y_m(s)\> -
{{2\eps_m}\over {\phi^3(x_m(s))}}{\|\nabla\phi(x_m(s))\|}^2\leq} & \nonumber \\
& {\displaystyle
H_\phi(x_m(s))[\dot x_m(s),\dot x_m(s)] -
H_\phi(y_m(s))[\dot y_m(s),\dot y_m(s)] + }& \nonumber \\
& {\displaystyle
M\phi(y_m(s))\<\dot y_m(s),\dot y_m(s)\> }&
\end{eqnarray}
where $y_m$ is defined by
(\ref{2.17}) and (\ref{2.18}).
In the following we shall always assume $s\in J$ and
$m\in\N$ large enough and we shall denote
by $C_1,\ldots ,C_{11}$ suitable positive constants.
Now by (\ref{3.6})
\begin{eqnarray} \label{3.20}
& {\displaystyle
H_\phi(x_m(s))[\dot x_m(s),\dot x_m(s)] -
H_\phi(y_m(s))[\dot y_m(s),\dot y_m(s)] = } & \nonumber \\
& {\displaystyle
L(x_m(s))[\dot x_m(s),\dot x_m(s)] -
L(y_m(s))[\dot y_m(s),\dot y_m(s)] =} & \nonumber \\
& {\displaystyle
L(x_m(s))[\dot x_m(s),\dot x_m(s)] -
L(y_m(s))[\dot x_m(s),\dot x_m(s)]  +}& \nonumber \\
& {\displaystyle
L(y_m(s))[\dot x_m(s) + \dot y_m(s),\dot x_m(s) - \dot y_m(s)].} &
\end{eqnarray}
Since $\phi\in{\cal C}^3$, $\eta\in{\cal C}^2$, using the
mean value theorem, the boundness of
$\(\|x_m \|_\infty\)_{m\in \N}$ and $\(\|y_m \|_\infty\)_{m\in \N}$ and {\it (iii)} of Theorem \ref{t1}, there results
\begin{eqnarray}  \label{3.11}
& {\displaystyle
L(x_m(s))[\dot x_m(s),\dot x_m(s)] -
L(y_m(s))[\dot x_m(s),\dot x_m(s)] \leq} & \nonumber \\
& {\displaystyle
C_1{|\dot\xm|}^2u_m(s).} &
\end{eqnarray}
Since $\eta_x(0,x)[v] = v$,
\begin{eqnarray}  \label{2.24}
& {\displaystyle
\dot x_m(s) - \dot y_m(s) =
\eta_x(0,x_m(s))[\dot x_m(s)] - } & \nonumber \\
& {\displaystyle
\eta_x(\tau_m(s), x_m(s))[\dot x_m(s)] -
\eta_s(\tau_m(s), x_m(s))\dot u_m(s),} &
\end{eqnarray}
so by (\ref{2.24}) the mean value theorem and {\it (iii)} of Theorem \ref{t1}
\begin{eqnarray}\label{3.12}
& {\displaystyle
L(y_m(s))[\dot x_m(s) + \dot y_m(s),\dot x_m(s) - \dot y_m(s)]=} & \nonumber \\
& {\displaystyle
L(y_m(s))[\dot x_m(s) + \dot y_m(s),
\eta_x(0,x_m(s))[\dot x_m(s)] -
\eta_x(\tau_m(s), x_m(s))[\dot x_m(s)]] -} & \nonumber \\
& {\displaystyle
L(y_m(s))[\dot x_m(s) + \dot y_m(s),
\eta_s(\tau_m(s), x_m(s))\dot u_m(s)]\leq} & \nonumber \\
& {\displaystyle
C_2{|\dot\xm|}^2u_m(s) + B_m(s)\dot u_m(s)} &
\end{eqnarray}
where
\be\label{3.13}
B_m(s) =
- L(y_m(s))[\dot x_m(s) + \dot y_m(s),\eta_s(\tau_m(s),x_m(s))]
\ee
is a ${\cal C}^1$ function.
By (\ref{2.11}), (\ref{3.9}), Lemma \ref{l3}  and
{\it  (ii)} of Theorem \ref{t1} we get
\be\label{3.14}
\<\dot \xm,\dot\xm\>\leq C_3
\ee
for a suitable constant $C_3$, hence by (\ref{3.10}),
(\ref{3.11}), (\ref{3.12}) and (\ref{2.27}) we get
\be \label{3.15}
\ddot u_m(s)\leq C_4 u_m(s) + B_m(s)\dot u_m(s).
\ee
Note that by (\ref{3.14}) and {\it (iii)} of Theorem \ref{t1}
\be\label{3.16}
|B_m(s)|\leq C_5
\ee
and by (\ref{2.5}), {\it  (ii)} of Theorem \ref{t1}
$$
|D_s\dot\xm|\leq C_6.
$$
Since
\begin{eqnarray*}
& \ddot\xm = P(\xm)\ddot\xm + Q(\xm)\ddot\xm &   \\
& =D_s \dot \xm -dQ ( \xm )[ \dot \xm , \dot \xm ] &
\end{eqnarray*}
where $P,Q$ are as in (\ref{3.2}), (\ref{3.3}), we also
get
$$
|\ddot x_m|\leq C_7.
$$
Then by using {\it (iii)} of Theorem \ref{t1}, standard arguments show that
\be\label{3.18}
|\dot B_m(s)|\leq C_8.
\ee
Therefore by (\ref{3.15}), (\ref{3.16}) and
(\ref{3.18})
integrating by parts, for $s>s_m$ there results
\begin{eqnarray}\label{3.19}
& {\displaystyle
\dot u_m(s) =
\integ\ddot u_m(\tau)d\tau\leq
C_4\integ u_m(\tau)d\tau +
\integ B_m(\tau)\dot u_m(\tau)d\tau =} & \nonumber \\
& {\displaystyle
C_4\integ u_m(\tau)d\tau + B_m(s)u_m(s) -
B_m(s_m)u_m(s_m) -
\integ \dot B_m(\tau)u_m(\tau)d\tau \leq}  & \nonumber \\
& {\displaystyle
C_9 \integ u_m(\tau)d\tau +
C_5 u_m(s) + C_5 u_m(s_m). } &
\end{eqnarray}
Hence
\begin{eqnarray}\label{3.21}
& {\displaystyle
u_m(s)\leq u_m(s_m) +
C_9 \integ\(\int_{s_m}^\tau u_m(r)dr\)d\tau +
C_5 \integ u_m(\tau)d\tau + C_5 u_m(s_m)\leq} & \nonumber \\
& {\displaystyle
C_{10} u_m(s_m) + C_{11}\integ u_m(\tau)d\tau.} &
\end{eqnarray}
By the Gronwall lemma we get
$$
u_m(s)\leq C_{10} u_m(s_m)\exp (C_{11}(s-s_m) )
$$
so by (\ref{2.14}) we have
$$
u_m\longrightarrow 0 \hbox{ uniformly}
$$
getting a contradiction. $\Box$
\vspace{0.5 cm}

\noindent
{\bf Proof of Theorem \ref{t1}.} By Remark \ref{r1} and Proposition \ref{l2}, we can find a family $\(x_\eps\)_{\eps > 0}$ of critical points
of $f_\eps$ such that (\ref{2.13}) holds.
By (\ref{2.12}) and the Poincar\'e inequality,
$\(x_\eps\)_{\eps > 0}$
is bounded in
$H^{1,2} ([0,1], \R^N )$, so there is a subsequence
$\(x_{\eps_m}\)_{m\in \N}$ such that
\be\label{2.29}
x_{\eps_m}\longrightarrow x
 \hbox{ weakly in }\,
  H^{1,2} ([0,1], \R^N ),
\ee
where $x\in\Om\(\D\)$ since the convergence is also uniform and (\ref{2.13}) holds. Now it is easy to prove that $f$ attains a minimum value at $x$. Indeed by
(\ref{2.29}) and (\ref{feps})
$$
f(x)\leq
\liminf_{m\longrightarrow\infty} f_{\eps_m}(x_{\eps_m})\leq f(y)
$$
for any $y\in\Om\(\D\)$, so $\D$ is convex.
Finally, if $\D$ is not contractible in itself, the proof
can be carried out exactly as the one of Theorem 0.2 of \cite{sal1}. $\Box$
\vspace{0.5 cm}

\noindent
{\bf Proof of Theorem \ref{t2}.} The only difference with the proof of Theorem \ref{t1}
concerns the a priori estimates of Proposition \ref{l2} that here are simpler
because they do not require a projection. Indeed, defining
$x_m$, $u_m$, $s_m$, $J$ as in Proposition \ref{l2}, we get by {\it  (ii)} and
 {\it (iii)} of Theorem \ref{t2}, for $m$ sufficiently large
$$\ddot u_m (s) \leq M \phi (x_m (s))\< \dot x_m (s), \dot x_m (s) \>
- \frac{2 \eps_m a^2}{\phi^3 (x_m (s))}.$$
As (\ref{2.11}) holds, we get
$$\ddot u_m (s) \leq M_1 u_m (s) + M_2 \frac{\eps_m}{u_m(s)}
-M_3 \frac{\eps_m}{u_m^3(s)}$$
for $M_1 , M_2 , M_3 >0$, so that, if $\delta$ is sufficiently small
$$\ddot u_m (s) \leq M_1 u_m (s).$$
Then, by the Gronwall Lemma, we immediately get a contradiction. $\Box$

\section{Applications}
\label{Appendix}

\noindent {\bf (a) Some examples.}
(1) First, we will check that Theorems \ref{t1}, and \ref{t2} can be applied
in cases where neither the elementary considerations for differentiable boundary nor Theorem \ref{t0} are appliable.
Let $(\M,\<\cdot,\cdot\>)$ be a cylinder $C$ in $\R^3$ and let
$H$ be an helix in $C$. Set $\D = C\setminus H$ and take $\phi$ equal to the
distance to $H$ on a strip $S$ around $H$. Recall that $\D$ cannot be considered as a manifold with boundary $\partial \D = H.$
Clearly, Theorems \ref{t1} and \ref{t2} are applicable and, moreover,
we can perturb the metric  to make the constant $M$ in (\ref{1.1})
positive (for example, by conformally changing the metric
symmetrically around $H$ on $S$, see formula (\ref{jhes})); thus,
Theorem \ref{t0} may be not applicable now.

\vspace{0.5 cm}

(2) Nevertheless, the following example shows that Theorem \ref{t0}
may be applicable when Theorems \ref{t1}, \ref{t2} are not, even if the condition for the Hessian in Theorem \ref{t2} is satisfied.
Consider $\D = ]0,\infty[\times ]0,\infty[$, and $\M = \R^2$
equipped with the Euclidean metric (this metric can be perturbed as
suggested in last example to make this one
non--trivial; note that conditions {\it (ii)}, {\it (iii)} in Theorem \ref{t1} are independent of the metric on $\M$).
If $\phi(x,y) = \sqrt{xy}$, then
$${\|\nabla\phi(x,y)\|}^2 = {{x^2 + y^2}\over{4xy}}> {1\over 4},$$
so the first inequality in {\it (ii)} of Theorems \ref{t1}, \ref{t2} is satisfied, but not the second one. On the other hand, if we choose $\hat \phi (x,y) = xy$ then
$${\|\nabla \hat \phi(x,y)\|}^2 = x^2 + y^2,$$
which satisfies the upper local bound around each point,
but not the lower one (compare with Remark \ref{rem1.7}(2)). Nevertheless,
it is straightforward to check that Theorem \ref{t0} is applicable for $\phi$ as well as for $\hat \phi$.

\vspace{0.5 cm}

(3) It is clear that condition {\it (iii)} of Theorem \ref{t2} on all tangent vectors $v$  may not hold under the corresponding hypothesis {\it (iv)} of Theorem \ref{t1}.  The following example shows that Theorem \ref{t2} may be appliable when Theorem \ref{t1} is not, because of the hypothesis {\it (iii)} of this theorem and the requirement that $\D$ must be a subset of a complete manifold.
Consider $\D=\M = \R\setminus \{0\}\times \R$ endowed with the Riemannian metric given in polar coordinates $r,\theta$ by
$$
dx^2 = dr^2 + {1\over{r^2}}d\theta^2.
$$
Choosen $\phi(r,\theta) = r$, it follows $\nabla\phi = \partial_r$, so $\|\nabla\phi (r, \theta ) \| = 1$.
Standard calculations show that the (normalized) flow of $\nabla\phi$ is
$\eta(s,(r,\theta)) = (r-s,\theta)$ and the norm of the partial
derivative $\eta_{(r,\theta)}$ is not bounded, so we cannot apply
Theorem \ref{t1}. Moreover, note that $\M$ is not complete and its curvature
along  incomplete radial geodesics diverges, so $\D$ is not isometric
to a domain of a complete Riemannian manifold ($\partial_c \D$
is topologically a circumference). However, it is easy to check that Theorem \ref{t2} is applicable.

\vspace*{0,5 cm}
\noindent {\bf (b) Trajectories of Lagrangian systems.}
The interest in the study of the geodesic connectedness of a complete Riemannian manifold is related to the existence of trajectories of a Lagrangian system joining two fixed points. More precisely, consider a potential $V\in C^2 (\M , \R )$ bounded from above in a domain $\D$ and the system
\be \label{lag}
\left\{ \begin{array}{ll}
D_s \dot x = -\nabla V(x) \\
x(0) = p,
  x (1) =  q.
\end{array}
\right.
\ee
Each solution $x: [0,1] \rightarrow \D$ of (\ref{lag}) has constant energy, that is, there exists $E\in\R$ such that, for any $s\in [0,1]$
$$\frac{1}{2} \< \dot x , \dot x \> + V(x) = E.$$
We fix $E> \sup_\D V$ and consider the Jacobi metric
\be \label{jac}
\< \cdot , \cdot \>_E = (E - V(x))\< \cdot , \cdot \> \quad \quad \hbox{on }  \overline \D.
\ee
As $\overline \D$ is complete for $\<\cdot , \cdot \>$ then it is also complete for  $\< \cdot , \cdot \>_E$ and, thus, we can extend $\< \cdot , \cdot \>_E$ to a complete Riemannian metric on all $\M$.
It is well--known that the geodesics on $\D$ with respect to the Riemannian
metric $\< \cdot , \cdot \>_E$ are, up to reparametrizations, solutions of (\ref{lag})
with energy $E$ and vice--versa. Let $\phi$ be a function satisfying conditions {\it (i), (ii), (iii)} in Theorem \ref{t1}, which are independent of the metric. If assumption {\it (iv)} is satisfied with respect to $\< \cdot , \cdot \>_E$ then the existence of at least one solution of (\ref{lag}) with energy $E$ is proved. Notice that the Hessian of $\phi$ with respect the two metrics is linked by  the following relation
\be\label{jhes}
H^E_\phi(x)[v,v]=  H_\phi(x)[v,v] + \<\nabla\phi(x),\nabla u(x)\>\<v,v\>
- 2\<\nabla u(x),v\>\<\nabla\phi(x),v\>
\ee
for any $x\in \M$, $v\in T_x \M$ where
$$u(x)=\frac{1}{2} \log (E-V(x)).$$
Thus, if {\it (iv)}  of Theorem \ref{t1} is verified with respect to $\< \cdot , \cdot \>$, then it is satisfied by $\< \cdot , \cdot \>_E$ when for each $y \in \partial \D$ there exists a neighborhood $U$ and a (positive) constant $M'\in \R$ such that
\be\label{rep}
\< \nabla \phi (x) , \nabla V(x) \> \geq -M' \phi(x) 
\ee
$\forall x\in \phi^{-1} (a_m) \cap U$. Summing up the following result holds.
\begin{cor}
Let $\(\M,\<\cdot,\cdot\>\)$ be a complete Riemannian manifold and
assume that there exists a positive and differentiable function $\phi$ on $\D$ satisfying {\it (i)--(iv)} of Theorem \ref{t1}. Then, if $V\in C^2 (\M , \R )$ is bounded  from above on $\D$ and (\ref{rep}) holds, for any
$E>\sup_\D V$ there exists a solution with energy $E$ of (\ref{lag}).
Moreover, if $\D$ is non contractible in itself, then  for any
$E>\sup V$ there exist infinitely many solutions of (\ref{lag}) with energy $E$.
\end{cor}

\vspace{0.5 cm}
\noindent{\bf Acknowledgment.} We wish to thank Domingo Rodr\'{\i}guez for his kind help for the figures of this paper.


\begin{figure} [ht]
\vskip 7.5cm \includegraphics{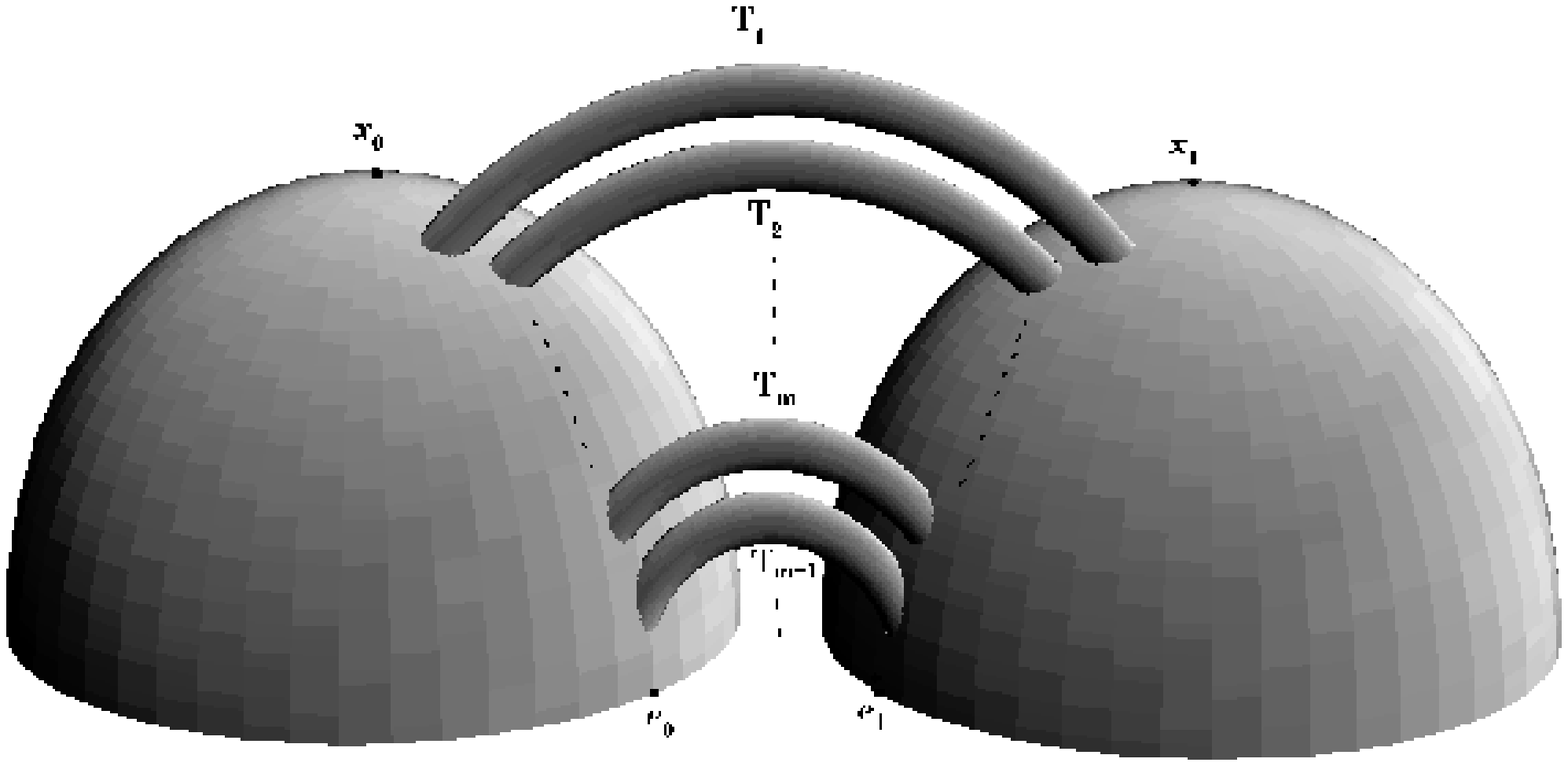} \label{Fig1} \caption{The two hemispheres $H_0, H_1$
are connected by a sequence of inmersed tubes $(T_m)$ such
 that the length of a minimizing connecting curve $\gamma_m$  through $T_m$ is bigger
 than the length of the corresponding $\gamma_{m+1}$ through $T_{m+1}$.}
\end{figure}

\begin{figure} [ht]
\vskip 7.5cm \includegraphics{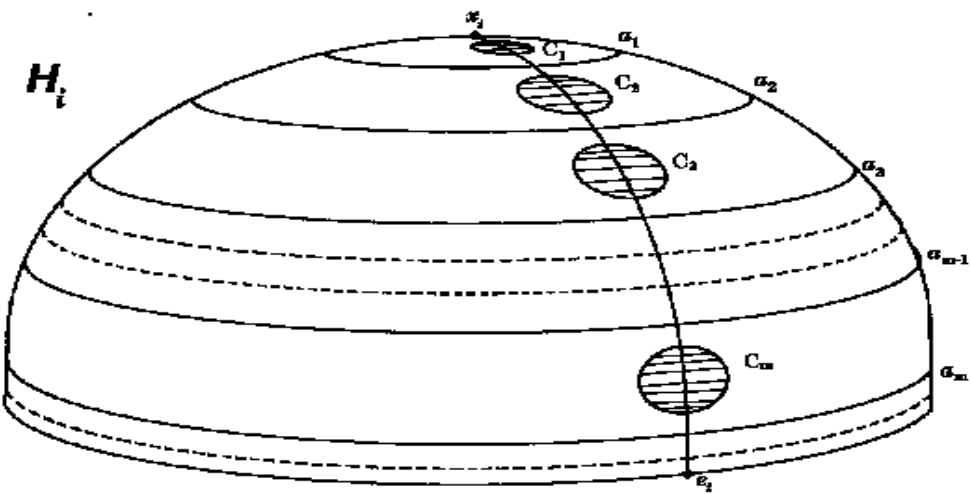} \label{Fig2} \caption{The dashed circles $C_m$ are
removed to attach the tubes $T_m$. The centers of the mouths of
the tubes lie in a meridian $\Upsilon$ and converge to an
equatorial point $e_i$ ($i=0,1$). Each parallel with heigth $a_m$
is equal to the connected component of the boundary $\partial
\D_m\ \cap H_i$, and it is clearly convex.}
\end{figure}

\end{document}